\documentstyle{article}
\newcommand{\bd}{\begin{defn}\em}
\newcommand{\bdl}[1]{\begin{defn}\em\label{#1}}
\newcommand{\bl}{\begin{lem}\em}
\newcommand{\bll}[1]{\begin{lem}\em\label{#1}}
\newcommand{\bp}{\begin{prop}\em}
\newcommand{\bpl}[1]{\begin{prop}\em\label{#1}}
\newcommand{\bt}{\begin{thm}\em}
\newcommand{\btl}[1]{\begin{thm}\em\label{#1}}
\newcommand{\bel}[1]{\begin{ex}\em\label{#1}}
\newcommand{\bc}{\begin{cor}\em}
\newcommand{\bcl}[1]{\begin{cor}\em\label{#1}}
\newcommand{\bproof}{\begin{pruuf}}
\newcommand{\eproof}{\end{pruuf}}
\newcommand{\ed}{\end{defn}}
\newcommand{\el}{\end{lem}}
\newcommand{\ep}{\end{prop}}
\newcommand{\et}{\end{thm}}
\newcommand{\ee}{\end{ex}}
\newcommand{\ec}{\end{cor}}

\newcommand{\ben}{\begin{enumerate}}
\newcommand{\een}{\end{enumerate}}
\newcommand{\bit}{\begin{itemize}}
\newcommand{\eit}{\end{itemize}}
\newcommand{\bds}{\begin{description}}
\newcommand{\eds}{\end{description}}

\newcommand{\crl}{\mbox {$\,\ -\kern-0.5em <\ $}}

\newcommand{\delete}[1]{} 
\newcommand{\suppress}[1]{} 
\newcommand{\excl}[1]{} 

\newtheorem{thm}{Theorem}
\newtheorem{lem}{Lemma}
\newtheorem{cor}{Corollary}
\newtheorem{ex}{Example}
\newtheorem{prop}{Proposition}
\newtheorem{defn}{Definition}

\newenvironment{pruuf}{\noindent{\bf Proof }}{}

\newcommand{\doar}{\mbox {$\,\ \rightarrow\kern-0.7em\rightarrow\ $}}

\newcommand{\calR}{{\cal R}}

\begin{document}
\title{Limiting values of rational functions at 
the points of discontinuity}

\author{Yaacov Tzeitlin  \\
}

\date{16.02.2000}
\maketitle


\section{Introduction and conventions}

In the present article we describe a class of algebraic curves on which
rational functions of two arguments may reach all their possible 
limiting values. We also
solve a similar question for functions that can be represented as a uniform
limit of a sequence of bounded rational functions.
\vskip 0.3cm
Some standart notation. Below we denote by 
$\mathbf{N}, \mathbf{Z^+}, \mathbf{Q}, \mathbf{R}$
the sets of all natural, nonnegative integer, 
rational and real numbers, respectively.
Further, all functions are real valued.
 
\begin{defn} \label{def1}
Let $f(x,y)$ be a real valued function of two real arguments. 
We say that a point $(x_0,y_0)$ is a point of 
bounded discontinuity if $f$ is defined and bounded 
in a certain nbd of a point $(x_0,y_0)$ and $(x_0,y_0)$ is a point of
the discontinuity. 
\end{defn}

Let $f(x,y)=\frac{p(x,y)}{q(x,y)}$ be a rational function when 
$p(x,y), q(x,y)$ are polinomials. 
For every $(x_0, y_0)$ there exists the limit
$$\lim_{x\to x_0} \lim_{y\to y_0} f(x,y).$$ 
If
$\lim_{x\to x_0} \lim_{y\to y_0} f(x,y) \ne \infty$, then we will suppose
that $f(x_0, y_0)=\lim_{x\to x_0} \lim_{y\to y_0} f(x,y)$. That is, 
we assume that such 
rational functions are defined at the points of bounded discontinuity.
Considering the rational function $\frac{p(x,y)}{q(x,y)}$ we
will asume below that $p(x,y)$ and $q(x,y)$ are mutually prime.
Let us introduce some additional notation: 
\newline
$\calR =$ \{the set of all rational functions with two arguments\}
\newline $\calR_{(x,y)}=\{ f(x,y)\in \calR | f(x,y) $
is bounded in some nbd of $(x,y)$\}.
\newline $\calR_G=\bigcap \{\calR_{(x,y)} | (x,y)\in G \}$ 
for a certain
$G\subset \mathbf{R}^2$.

Note that $\calR_{(x,y)}$ and $\calR_G$ are commutative real algebras with 
respect to pointwise operations. 

\begin{lem} \label{lem1}
If $f(x,y)\in \calR$ then the set of all points of 
bounded discontinuity is finite. 
\end{lem}
  

\bproof
Indeed if $f(x,y)=\frac{p(x,y)}{q(x,y}$, then the set of all bounded 
discontinuity is a subset of all roots of the following system
$$p(x,y)=0, \quad q(x,y)=0.$$

On the other hand the number of such solutions is at most
$deg(p) \cdot deg(q)$, where $p$ and $q$ are mutually prime.
\hfill $\Box$
\eproof

It is important to know what is the behaviour of the function $f(x,y)$ 
near the point $(x_0, y_0)$. This can be described by the set 
$D_f(x_0,y_0)$ of all limiting values of $f$ at $(x_0, y_0)$.
By the definition,
$$D_f(x_0,y_0)= \bigcap \{\overline{f(V(x_0,y_0)} | V\in N_{(x_0, y_0)} \},$$
 where $N_{(x_0, y_0)}$ denotes the set of all nbd's of the point 
$(x,y)$.
Let $\underline{f}=inf D_f(x_0, y_0), \overline{f}=sup D_f(x_0, y_0)$.

\begin{lem} \label{lem2}
If $f(x,y)\in \calR_{(x_0, y_0)}$, then $D_f^{(x_0, y_0)}
=[\underline{f}, \overline{f}]$.
\end{lem}

\bproof 
Observe that 
$$\underline{f}\in D_f^{(x_0, y_0)},  \quad 
\overline{f} \in D_f^{(x_0, y_0)}$$
Let $d\in (\underline{f}, \overline{f})$. Then by Lemma \ref{lem1},
$f$ is continuous at every point $\not= (x_0, y_0)$ of 
a sufficiently small nbd of $(x_0, y_0)$. Therefore, 
by \emph{the intermediate 
value theorem} there exists a point $(x_d, y_d)$ such that $f(x_d, y_d)=d$.
\hfill $\Box$
\eproof

\begin{lem} \label{lem3}
Let $G$ be a domain in $\mathbf{R}^2$ and 
$f(x,y)=\frac{p(x,y)}{q(x,y)} \in \calR_G$. Then $q$ has a constant sign 
on $G$.
\end{lem} 

\bproof
Since $G$ is an open path-connected subset of $\mathbf{R}^2$, 
the same is true for the 
subset $\{(x,y) \in G | q(x,y) \not=0\}$ and arbitrary pair 
$(x_1, y_1), (x_2, y_2)$ of points from this subset can 
be joined in this 
open subspace by a path. This implies that $q(x_1, y_1), q(x_2, y_2)$ 
have the same sign.
\hfill $\Box$
\eproof

Further we will deal with the characters on the algebra $\calR_{(x_0,y_0)}$.
Recall that a (nontrivial) real \emph{character} on an algebra $A$ 
is a (nontrivial) homomorphism from $A$ into $\mathbf{R}$ 
(see \cite{[2]}). Denote by $X^A$ the set of all nontrivial characters.
Endow it with the \emph{weak topology} $\sigma (X^A, A)$, that is, 
topology generated by all maps $\tilde{a}: X^A \to \mathbf{R}, 
\tilde{a}(\chi)=\chi(a)$.

We need the following
\bd \label{def2}  

1. Let $\{(x_{\gamma}, y_{\gamma})\}_{\gamma \in \Gamma}$ be a net \cite{[3]}
that converges to the point $(x_0, y_0)$ in $\mathbf{R}^2$ such that 
$$\forall f\in \calR_{(x_0, y_0)} \quad \exists lim_{\gamma \in \Gamma} 
f(x_{\gamma}, y_{\gamma}) \not=\infty.$$ 
 
The formula $\chi (t)= \lim_{\gamma \in \Gamma} f(x_{\gamma}, y_{\gamma})$
defines a character $\chi$ of the algebra $\calR_{(x_0,y_0)}$. 
We will say that the net $\{(x_{\gamma}, y_{\gamma})\}_{\gamma \in \Gamma}$ 
determines the character $\chi$.

2. Let $L$ be a curve in $\mathbf{R}^2$ and $(x_0,y_0)\in L$. then if

$$\forall f\in \calR_{x_0,y_0} \quad
\exists \lim_{(x,y)\in L \to (x_0,y_0)} f(x,y)$$
then the formula
$\chi (f)=\lim_{(x,y)\in L \to (x_0,y_0)} f(x,y)$ defines a character of the 
algebra $\calR_{x_0,y_0}$. We will say that the curve $L$ determines the 
character $\chi$.

3. Let $\{L_N\}_{N\in \mathbf{N}}$ be a sequence of curves such that each of 
them determines the character $\chi_{N}$. Then if 
$\{L_N\}_{N\in \mathbf{N}}$ converges in
$\sigma (\chi^{\calR_{(x_0,y_0)}}, \calR_{(x_0,y_0)})$ 
to the character $\chi$,
then we say that $\{L_N\}_{N\in \mathbf{N}}$ determines $\chi$.
\ed

In the following definition we introduce functions and series 
that allow us to describe the characters of the algebra $\calR_{(x_0,y_0)}$.

\bd \label{def3}

(i) For every $N\in \mathbf{Z^+}$ define the function $\Phi_N(t)$ by
$$\forall N\in {\bf N} \quad \Phi_N(t)=\sum_{k=1}^N \phi_k|t|^{e_k}$$
where $\phi_k \in {\bf R}$, $\phi_k \not= 0$, $e_k \in {\bf Q}$ and
$1\le e_1 < \cdots < e_N$. If $N=0$ then by definition $\Phi_0(t)=0$.
All such functions will be called functions of type (i).

(ii) We say that the series $\sum_{k=1}^{\infty} \phi_k|t|^{e_k}$
is of type (ii), if all partial sums are of type (i). 

\ed 

Denote by $\cal{I}$ the set of all functions of the type (i) and by 
$\cal{II}$ the set of all series of the type (ii). 
Finally, denote by $\cal{II\infty}$ the set of all series
$\{\sum_k \phi_k |t|^{e_k} \in {\cal II} 
| \lim_{k \to \infty} e_k =\infty \}$.

We naturally can define on the sets $\cal{I}, {\cal II \infty}$ 
the binary operations: addition 
and multiplication. Namely, 
let $\sum_i \alpha_i |t|^{m_i}$ and $\sum_j \beta_j |t|^{r_j}$ 
belong to $\cal{I}\cup \cal{II}\infty$. Suppose that
$$\sum_i \alpha_i |t|^{m_i} + \sum_j \beta_j |t|^{r_j}=
\sum_k \phi_k |t|^{e_k},$$
where $\{e_k\}$ is a sequence which is obtained by writing elements of
$m_i, r_j$ in their natural order and
$\phi_k=\left\{ \begin{array}{ll}
                \alpha_i & \mbox{if $e_k=m_i \notin \{r_j\}$} \\
                \beta_j & \mbox{if $e_k=r_j \notin \{m_i\}$} \\
                \alpha_i + \beta_j & \mbox{if $e_k=m_i=r_j$}    
                 \end{array}
                \right.  $             

The multiplication is defined by 

$ (\sum_i \alpha_i |t|^{m_i}) \cdot (\sum_j \beta |t|^{e_k})=
\sum_k \phi_k |t|^k.$ 
Here $\phi_k=\sum_{m_i+r_j=e_k} \alpha_i \beta_j $ and 
$\{e_k\}$ is obtained by rewriting the elements of 
$\{m_i + r_j\}_{i,j}$ as an increasing sequence.
Thus $\cal{I}\cup \cal{II\infty}$ is an algebra over reals.
It is impossible to extend these operations on 
$\cal{I}\cup {\cal II}$. 
However such extension exists for natural degrees.

\begin{lem} \label{lem4}
For every $\sum_k \phi_k |t|^{e_k} \in \cal{I}\cup \cal{II} $
and $n\in N$ is defined the following equality
$(\sum_k \phi_k |t|^{e_k})^n$ $=\sum_m g_m |t|^{s_m}$. 
Here, $s_m$ is obtained
by rewriting the elements of 
$\{e_{k_1}+\cdots +e_{k_n} |  k_1\leq k_2 \cdots \leq k_n\}$ 
as an increasing sequence and 
$g_m=\sum_{k_1+k_2+\cdots +k_n=s_m} \phi_{k_1}\cdots \phi_{k_n}$. 
\end{lem} 

\bproof 
For $\sum_k \phi_k |t|^{e_k} \in \cal{I}\cup \cal{II\infty} $
the assertion is true because $ \cal{I}\cup \cal{II\infty} $
is an algebra.
It suffices to prove for the series $\sum_k \phi_k |t|^{e_k} $
of the type $\cal{I}$ such that 
$\lim_{k\to \infty} e_k =e \neq \infty$.

We have to show that 
$$\sum_{e_{k_1}+ \cdots +e_{k_n}=s_m} \phi_{k_1}\cdots \phi_{k_n}$$ 
is well-defined for every 
$s\in \{e_{k_1}+\cdots+e_{k_n} | k_1\leq \cdots k_n\}$ 
and that 
$\{e_{k_1}+\cdots+e_{k_n} | k_1\leq \cdots k_n\}$
can be represented as an increasing sequence
$\{s_m\}_{m\in N}$ for each $n\in N$.

In order to show this it is sufficient to establish that the 
number of the roots of the equation 
$e_{k_1}+\cdots+e_{k_n}=s$
is finite for every
$s\in \{e_{k_1}+\cdots+e_{k_n} | k_1\leq \cdots \leq k_n\}$.

We use the induction on $n$. For $n=1$ the assertion is trivial.
We have to check the case of $n+1$ assuming that the assertion is true for 
$n$. Consider the equation
$e_{k_1}+\cdots+e_{k_n}+e_{k_{n+1}}=s$.
Then 
$(n+1)e_{k_1} < e_{k_1} + \cdots + e_{k_{n+1}} < (n+1)e_{k_{n+1}}<(n+1)e$
and $e_{k_1}<\frac{s}{n+1}<e$.
However, for sufficiently big $N$ and for $k \geq N$ holds 
$e_k>\frac{s}{n+1}$. Hence we obtain that 
$k_1<N$.
By inductive assumtion, the equation
$e_{k_2}+\cdots+e_{k_{n+1}}=s-e_{k_1}$
has finitely many roots for every $k_1<N$. Therefore,
the equation $s=e_{k_1}+\cdots+e_{k_n}+e_{k_{n+1}}$ 
also has finite number of the roots. This proves the lemma. \hfill $\Box$
\eproof

\vskip 0.3cm
Making use Lemma \ref{lem4}, 
we can obtain an assymptotic characterization of the behaviour of 
a polinomial on partial sums of the series of type  
${\cal II}$.

\begin{lem} \label{lem5}
Let $P(x,y)$ be a polinomial, 
$\sum_{k=1} ^{\infty} \phi_k |t|^{e_k}\in \cal{II}$
and let 
$\{\Phi_N(t)\}_{N\in \cal{N}}$ 
be the sequence of all partial sums of the above series. Then there exists
$N_0 \in \cal{N}$
such that for every $N>N_0$ we have  
$p(x,y_0+\Phi_N(x-x_0))=p(x_0,y_0)+c(x-x_0)^u+o((x-x_0)^u)$
for $x>x_0$ and
$p(x,y)=p(x_0,y_0)+d(x_0-x)^v+o((x_0-x)^v)$ 
for $x<x_0$. Here 
$u,v$ are rational and $c,d$, real numbers.  
\end{lem} 

\section{Characters of the algebra ${\cal R}_{(x_0,y_0)}$}

Now we are in the position to describe curves and sequences of curves that
determine characters of the algebra $\calR_{(x_0,y_0)}$.

\begin{thm} \label{Th1}
Let $\Phi_N(x-x_0)$ be an arbitrary function of type $\cal{I}$. Then the
curves and sequences of curves with the following equations 
determine characters of the algebra $\calR_{(x_0,y_0)}$.

\begin{enumerate}
\item  $y=y_0+\Phi_N(x-x_0)$ \hskip 0,5cm  $N\in {\bf Z}^+$. 

\item $y=y_0+\Phi_N(x-x_0)=|x-x_0|^e$, where $N\in {\bf Z}^+$ and 
$e\notin {\bf Q}$. We assume in addition that $1\leq e$ if $N=0$ and
$e_N \le e$ if $N\in {\bf N}$. 

\item 
$y=\left\{ \begin{array}{ll}
                y_0+\Phi_N(x-x_0)+\frac{|x-x_0|^e}{ln|x-x_0|}
& \mbox{if $x\neq x_0$} \\
                y_0 & \mbox{if $x=x_0$} \\
                                 \end{array}
                \right.  $

Here and in (4) we assume that $N\in {\bf Z}^+, e\in Q$, $1\leq e$
if $N=0$ and $e_N\le e$ for $N\in {\bf N}$.

\item             
$y=\left\{ \begin{array}{ll}
                y_0+\Phi_N (x-x_0)+|x-x_0|^e ln |x-x_0| 
& \mbox{if $x\neq x_0$} \\
                y_0 & \mbox{if $x=x_0$} \\
                                 \end{array}
                \right.  $             

\item The sequence 
$\{y_0+\Phi_N(x-x_0)\}_{N\in {\bf N}}$, where 
$\{\Phi_N(x-x_0)\}_{N\in {\bf N}}$ is the sequence of partial sums of series 
from ${\cal II}$.  

\end{enumerate}
Each of above equations from (1,2.3,4) and the series from (5) define, in 
general, different characters for the cases $x_0\leq x$ and $x\leq x_0$.  

\end{thm}
\bproof
In the cases (1,2,3,4) it is straightforward to show that for every
$f(x,y)\in \calR_{(x_0,y_0)}$ there exists a limit which determines 
the corresponding character. In the case (5), we use also 
Lemma \ref{lem5}.
\hfill $\Box$ 
\eproof

Analogously can be established

\begin{thm} \label{Th1'}
Let $\Phi_N(y-y_0)$ be an arbitrary function of type ${\cal I}$. Then the
curves and sequences of curves with the following equations 
determine characters of the algebra $\calR_{(x_0,y_0)}$.

\begin{enumerate}
\item  $x=x_0+\Phi_N(y-y_0)$ \hskip 0,5cm  $N\in {\bf Z}^+$. 

\item $x=x_0+\Phi_N(y-y_0)=|y-y_0|^e$, where $N\in {\bf Z}^+$ and $e$ be 
irrational. We assume in addition that $1\leq e$ if $N=0$, and
$e_N \le e$ if $N\in \cal{N}$. 

\item 
$x=\left\{ \begin{array}{ll}
                x_0+\Phi_N(y-y_0)+\frac{|y-y_0|^e}{ln|y-y_0|} 
& \mbox{if $y\neq y_0$} \\
                x_0 & \mbox{if $y=y_0$} \\
                                 \end{array}
                \right.  $

Here and in (4) we assume that $N\in {\bf Z}^+, e\in Q$. Moreover $1\leq e$
if $N=0$ and $e_N\le e$ for $N\in {\bf N}$.

\item             
$x=\left\{ \begin{array}{ll}
                y_0+\Phi_N(x-x_0)+|x-x_0|^e ln|x-x_0|
& \mbox{if $x\neq x_0$} \\
                y_0 & \mbox{if $x=x_0$} \\
                                 \end{array}
                \right.  $             

\item The sequence 
$\{x_0+\Phi_N(y-y_0)\}_{N\in {\bf N}}$, where 
$\{\Phi_N(y-y_0)\}_{N\in {\bf N}}$ is the sequence of partial sums of series 
from ${\cal II}$.  

\end{enumerate}
Each of above equations from (1,2.3,4) and the series from (5) define, in 
general, different characters for the cases $y_0\leq y$ and $y\leq y_0$.  

\end{thm}

Now we describe the nets that determine the characters of the algebra
$\calR_{(x_0,y_0)}$. First of all some useful remarks about such nets.
 
Every such net necessarily tends to 
$(x_0,y_0)$ in $\mathbf{R}^2$. We can ignore the trivial case of (finally)
stationary nets. If 
$\{(x_{\gamma},y_{\gamma})\}_{\gamma \in \Gamma}$ is not a finally stationary
net, then there exists a cofinal subnet
$\{(x_{\gamma},y_{\gamma})\}_{\gamma \in \Gamma_0}, 
\hskip 0,3cm \Gamma_0 \subset \Gamma$ such that one of the following 
conditions is fulfilled:

\begin{itemize}
\item [a)] 
$\lim_{\gamma \in \Gamma_0} \frac{y_{\gamma}-y_0}{x_{\gamma}-x_0}=\varphi\neq
\infty$
\item [b)]
$\lim_{\gamma \in \Gamma_0} \frac{x_{\gamma}-x_0}{y_{\gamma}-y_0}=\varphi\neq
\infty$.
\end{itemize}
 
Indeed, otherwise 
$\lim_{\gamma \in \Gamma} f(x_{\gamma},y_{\gamma})$
does not exist for the function
$$f(x,y)=\frac{(x-x_0)^2-(y-y_0)^2-(x-x_0)(y-y_0)}{(x-x_0)^2+(y-y_0)^2} \in
\calR_{(x_0,y_0)}.$$
Moreover we can assume that in both cases of (a), (b) the denomenators 
preserve the sign.

\begin{thm} \label{Th2} 
Let the net 
$\{(x_{\gamma},y_{\gamma})\}_{\gamma \in \Gamma}$ tends to 
$(x_0,y_0)$ in $\mathbf{R}^2$ and one of the following conditions 
hold:
\begin{itemize}
\item [a)]
$\lim_{\gamma \in \Gamma_0} \frac{y_{\gamma}-y_0}{x_{\gamma}-x_0}=\varphi\neq
\infty$, 
where $x_0< x_{\gamma}$ for every $\gamma\in \Gamma$. Or, 
$x_0> x_{\gamma}$ for every $\gamma\in \Gamma$;
\item [b)]
$\lim_{\gamma \in \Gamma_0} \frac{x_{\gamma}-x_0}{y_{\gamma}-y_0}=\varphi\neq
\infty$,
where $y_{\gamma}<y_0$ for every $\gamma\in \Gamma$. Or,
$y_0> y_{\gamma}$ for every $\gamma\in \Gamma$.
\end{itemize}
Then this net determines a character of the algebra $\calR_{(x_0,y_0)}$
iff one of the following conditions hold:

\begin{enumerate}
\item
\begin{itemize}

\item [a)] 
$\exists \Phi_N (t) \in {\cal I} \quad \forall e\in {\bf R} \quad
\lim_{\gamma \in \Gamma} 
\frac{y_{\gamma}-y_0-\Phi_N(x_{\gamma}-x_0)}{|x_\gamma-x_0|^e}=0.$

\item [b)] 
$\exists \Phi_N (t) \in {\cal I} \quad \forall e \in {\mathbf R} \quad
\lim_{\gamma\in \Gamma} 
\frac{x_{\gamma}-x_0-\Phi_N (y_{\gamma}-y_0)}{|y_{\gamma}-y_0|^e}=0.$

\end{itemize}
\item
\begin{itemize}

\item [a)]
$\exists \Phi_ N(t) \in {\cal I} \quad \exists e \notin {\mathbf Q}$
\newline
$\lim_{\gamma \in \Gamma} 
\frac{y_{\gamma}-y_0-\Phi_N(x_{\gamma}-x_0)}{|x_\gamma-x_0|^e}=0$ 
for $r<e$ 
and
\newline
$\lim_{\gamma \in \Gamma} 
\frac{y_{\gamma}-y_0-\Phi_N(x_{\gamma}-x_0)}{|x_\gamma-x_0|^e}=\infty$ 
for $r>e$. 

\item [b)] 
$\exists \Phi_N(t) \in {\mathcal I} \quad \exists e\notin \mathbf{Q}$
\newline
$\lim_{\gamma \in \Gamma} 
\frac{x_{\gamma}-x_0-\Phi_N(y_{\gamma}-y_0)}{|y_\gamma-y_0|^e}=0$ 
for $r<e$ and
\newline
$\lim_{\gamma \in \Gamma} 
\frac{x_{\gamma}-x_0-\Phi_N(y_{\gamma}-y_0)}{|y_\gamma-y_0|^e}=\infty$ 
for $r>e$. 
\end{itemize}

\item
\begin{itemize}
\item [a)]
$\exists \Phi_N(t) \in {\cal I} \quad \exists e\in \mathbf{Q}$
\newline
$\lim_{\gamma \in \Gamma} 
\frac{y_{\gamma}-y_0-\Phi_N(x_{\gamma}-x_0)}{|x_\gamma-x_0|^r}=0$ 
for $r\leq e$ 
and
\newline
$\lim_{\gamma \in \Gamma} 
\frac{y_{\gamma}-y_0-\Phi_N(x_{\gamma}-x_0)}{|x_\gamma-x_0|^r}=\infty$ 
for $r>e$. 

\item [b)] 
$\exists \Phi_N(t) \in {\cal I} \quad \exists e\in \mathbf{Q}$
\newline
$\lim_{\gamma \in \Gamma} 
\frac{x_{\gamma}-x_0-\Phi_N(y_{\gamma}-y_0)}{|y_{\gamma}-y_0|^r}=0$ 
for $r\leq e$ 
and
\newline
$\lim_{\gamma \in \Gamma} 
\frac{x_{\gamma}-x_0-\Phi_N(y_{\gamma}-y_0)}{|y_{\gamma}-y_0|^r}=\infty$ 
for $r>e$. 
\end{itemize}

\item
\begin{itemize}
\item [a)]
$\exists \Phi_N(t) \in {\cal I} \quad \exists e\in \mathbf{Q}$
\newline
$\lim_{\gamma \in \Gamma} 
\frac{y_{\gamma}-y_0-\Phi_N(x_{\gamma}-x_0)}{|x_\gamma-x_0|^r}=0$ 
for $r<e$ 
and
\newline
$\lim_{\gamma \in \Gamma} 
\frac{y_{\gamma}-y_0-\Phi_N(x_{\gamma}-x_0)}{|x_\gamma-x_0|^r}=\infty$ 
for $e\leq r$. 

\item [b)] 
$\exists \Phi_N(t) \in {\cal I} \quad \exists e\in \mathbf{Q}$
\newline
$\lim_{\gamma \in \Gamma} 
\frac{x_{\gamma}-x_0-\Phi_N(y_{\gamma}-y_0)}{|y_{\gamma}-y_0|^r}=0$ 
for $r< e$ 
and
\newline
$\lim_{\gamma \in \Gamma} 
\frac{x_{\gamma}-x_0-\Phi_N(y_{\gamma}-y_0)}{|y_{\gamma}-y_0|^r}=\infty$ 
for $e\leq r$. 
\end{itemize}

\item
\begin{itemize}
\item [a)]
$\exists \sum_{k=1}^{\infty} \phi_k |t|^{e_k} \in {\cal II} \quad
\forall N \in {\bf Z^+} \quad 
\lim_{\gamma \in \Gamma} 
\frac{y_{\gamma} -y_0 -\Phi_N (x_{\gamma} -x_0)}{|x_{\gamma} -x_0|^{e_N+1}}
=\phi_{N+1}$.
\newline
Here $\Phi_N (x_{\gamma} -x_0)=
\sum_{k=1}^{\infty} \phi_k |x_{\gamma} -x_0|^{e_k}$ 
for $N\in {\bf N}$
and $\Phi_0(x_{\gamma}-x_0)=0$.
\item [b)]
$\exists \sum_{k=1}^{\infty} \phi_k |t|^{e_k} \in {\cal II} 
\quad \forall N\in {\bf Z^+} \quad 
\lim_{\gamma \in \Gamma} 
\frac{x_{\gamma} -x_0 -\phi_N (y_{\gamma} -y_0)}
{|y_{\gamma} -y_0|^{e_{N+1}}}=\phi_{N+1}$. 
\newline
Analogously, here $\Phi_N (y_{\gamma} -y_0)=
\sum_{k=1}^{\infty} \phi_k |y_{\gamma} -y_0|^{e_k}$ 
for $N \in {\bf N}$
and $\Phi_0(y_{\gamma}-y_0)=0$.

\end{itemize}
\end{enumerate}
\end{thm}
\bproof
For simplicity assume 
$(x_0,y_0)=(0,0)$. We consider only the case of (a) with the assumption
$x_0\leq x_{\gamma}$. 
Other cases are very similar. 

We will present polinomials in the form 
$$p(x,y)=\sum_{k=0}^{n_p} (y-\Phi_N(x))^k x^{r_k}p_k(x),$$
where 
$x^{r_k}p_k(x)=\frac{1}{k!} \frac{\partial p}{\partial y^k}|_{y=\Phi_N(x)}=
x^{r_k}p_k(0)+0(x^{r_k}).$
Therefore 
$f(x,y)=\frac{p(x,y)}{q(x,y)} \in {\cal R} (x_0,y_0)$ we can
represent in the form
$$f(x,y)=\frac{\sum_{k=0}^{n_p} (y-\Phi_N(x))^kx^{r_k}p_k(x)}
{\sum_{k=0}^{n_q} (y-\Phi_N(x))^kx^{s_k}q_k(x)}.$$
Now we consider our conditions (1-5). If (1) holds, that is. if
$y_{\gamma}=0(x_{\gamma}^e)$ for $e\in {\bf R}$, then
$$\lim_{\gamma\in \Gamma} f(x_{\gamma},y_{\gamma})=
\lim_{\gamma\in \Gamma} 
\frac{\sum_{k=0}^{n_p} (y_{\gamma}-\Phi_N(x_{\gamma})^kx_{\gamma}^{r_k}
p_k(x_{\gamma})}{\sum_{k=0}^{n_q} (y_{\gamma}-\Phi_N(x_{\gamma})^k)
x_{\gamma}^{s_k}q_k(x_{\gamma})}
=\lim_{\gamma\in \Gamma} \frac{(y_{\gamma}-\Phi_N(x_{\gamma}))^n x_{\gamma}
^{r_n}p_n(x_{\gamma})}
{(y_{\gamma}-\Phi_N(x_{\gamma}))^m x_{\gamma}
^{s_m}q_m(x_{\gamma})}$$ 
where $n=min\{k| \frac{\partial ^k p}{\partial y^k} |_{y=\Phi_N(x)} \neq 0$
and  $m=min\{k| \frac{\partial ^k q}{\partial y^k} |_{y=\Phi_N(x)} \neq 0$.
Clearly, $\lim_{\gamma} f(x_{\gamma}, y_{\gamma}) \neq {\infty}$ exists
for every $f(x,y) \in {\cal R}_{(0,0)}$.
In order to examine the conditions (2), (3), (4), write $f(x,y)$ in the
form 
$$f(x,y)=\frac
{\sum_{k=0}^{m_p} (\frac{y-\Phi_N(x)}{x^e})^k x^{r_k+ke}p_k(x)}
{\sum_{k=0}^{n_q} (\frac{y-\Phi_N(x)}{x^e})^k x^{s_k+ke}q_k(x)}  
\quad \hskip 4cm (*)$$ 
Then 
$$\lim_{\gamma \in \Gamma} f(x_{\gamma},y_{\gamma})=
\lim_{\gamma \in \Gamma} \frac
{(\frac{y_{\gamma}-\Phi_N(x_{\gamma})}{x_{\gamma}^e})^n 
x_{\gamma}^{r_n+ne}p_n(x_{\gamma})}
{(\frac{y_{\gamma}-\Phi_N(x_{\gamma})}{x_{\gamma}^e})^m x_{\gamma}
^{s_m+me}}=
\left\{ \begin{array}{ll}
                \frac{p_n(0)}{q_n(0)}
& \mbox{if $m=n, r_n=s_n$} \\
                0 & \mbox{in other cases of (2)} \\
                                 \end{array}
                \right.   \hskip 1cm  (**)$$             

Now we give formulas like $(**)$ and methods for finding suitable 
$n,m$ for other cases.

For (2). Then $n$ and $m$ are determined by the conditions:
$$r_n+ne=min\{r_k+ke |\frac{\partial ^k p}{\partial y^k} |_{y=\Phi_N(x)}
\neq 0 $$

For (3). Then we set:
$n=min\{k| r_k+ke = u\}$ , where 
$u=min\{r_k+ke |\frac{\partial ^k p}{\partial y^k} |_{y=\Phi_N(x)} \neq 0 \}$
and $m=min\{k |s_k+ke=u\}$, where 
$u=min\{s_k+ke |\frac{\partial ^k q}{\partial y^k} |_{y=\Phi_N(x)} \}$.

For (4). Set:
$n=max\{k | r_k+ke =u\}$, where
$u=min\{r_k+ke | \frac{\partial ^kp}{\partial y^k}|_{y=\Phi_N(x)} \neq 0\}$,
$m=max\{k|s_k+ke=u\}$, where
$u=min\{s_k+ke|\frac{\partial ^kq}{\partial y^k}|_{y=\Phi_N(x)}\} \neq 0$.

Finally consider the case (5). Note that rewriting polinomials in the form
$$p(x,y)=\sum_{k=0}^{n_p} (y-\Phi_N(x))^k x^{r_k}p_k(x),$$
by Lemma \ref{lem5} 
we may state that starting from certain $N$ the value of
$r_k$ are unchanged for $k=0,1,\cdots,n_p$.
We will show next that for every polinomial $p(x,y)$ there exists a 
$k(p)\in {\cal N}$ such that for sufficiently big $N$ holds
$$p(x,y)=(\frac{y-\Phi_N(x)}{x^{e_{N+1}}})^{k(p)} x^{r_k(p)+k(p)e_{N+1}}+
0(x^{r_k(p)+k(p)e_N+1}) \quad \quad (***)$$
Indeed, if 
$\lim_{N\to {\infty}}e_N=\infty$ then clearly (***) 
holds for sufficiently big
$N$ and for 
$k(p)=min\{k| \frac{\partial ^k p}{\partial y^k} \neq 0\}$. Alternatively if 
$\lim_{N\to {\infty}}e_N=e <\infty$
and
$r=min\{r_k+ke| \frac{\partial ^k p}{\partial y^k} \neq 0\}$.
Then if 
$k(p)=max\{k|r_k+ke=r\}$ and $N$ is so big that 
$r<r_k+ke$ implies that
$r_{k(p)}+e_{N+1}k(p)<r_k+e_{N+1}k$,
then
$r_{k(p)}+e_{N+1}k(p)<r_k+e_{N+1}k$ whenever $k\neq k(p)$ and hence (***) 
holds. 
In order to finish the subcase of (5) observe that applying (***) to
$p(x,y),q(x,y)$ for
$f(x,y)=\frac{p(x,y)}{q(x,y)}\in {\cal R}_{(0,0)}$ the limit
$\lim_{\gamma \in \Gamma}f(x_{\gamma}, y_{\gamma})$ exists.
Therefore we have proved that each net which satisfies any of the
conditions (1-5) determines a character of the algebra 
${\cal R}_{(x_0,y_0)}$.
The present proof concerns all possibilities except the cases when
$e\in \mathbf{Q}$ and the limit (finite or infinite)
$\lim_{\gamma \in \Gamma} \frac{y_{\gamma}-\Phi_N(x_{\gamma})}
{x_{\gamma}^e}$ does not exist. But for these exceptional cases easily
can be constructed appropriate examples of
$f(x,y)\in {\cal R}_{(x_0,y_0)}$ such that 
$\lim_{\gamma \in \Gamma}f(x,y)$ does not exist. In particular,
if $\lim_{\gamma \in \Gamma} \frac{y_{\gamma}-\Phi_N(x_{\gamma})}
{x_{\gamma}^r}$
 exists for $r\neq e$
and  
if $\lim_{\gamma \in \Gamma} \frac{y_{\gamma}-\Phi_N(x_{\gamma})}
{x_{\gamma}^r}$ exists for 
$e=\frac{n_0}{m_0}$, then 
 $\lim_{\gamma \in \Gamma} f(x_{\gamma},y_{\gamma})$ does not exist for
$f(x,y)=\frac{y^{m_0}x^{n_0}}{y^{2m_0}+x^{2n_0}}$.
The proof is finished.
\hfill $\Box$
\eproof

\begin{cor}
If a net satisfies one the conditions i.(a) or i.(b) 
($i=1,2,3,4$) of Theorem \ref{Th2}, then it  
determines 
the same character as the curves described in i.
 of Theorem 1 or in i. of Theorem \ref{Th1'}. 
If a net satisfies (5), then the character 
determined by 
this net and the character determined by the sequence of the curves from  
Theorem \ref{Th1} or in Theorem \ref{Th1'} are the same. 
\end{cor}

\begin{cor}
If a character is determined by some 
net then such character can be determined
also by a {\it sequence}.
\end{cor}

\begin{cor}
\begin{enumerate}
\item
Let $\{(x_n,y_n)\}_{n\in {\cal N}}$ tends to $(x_0,y_0)$ then there exists 
a subsequence which determines a character of the algebra 
${\cal R}_{(x_0,y_0)}$.
\item
$\forall f(x,y)\in {\cal R}_{(x_0,y_0)} \quad 
\forall d\in D_{(x_0,y_0)}^f \quad
\exists \chi \in X^{{\cal R}_{(x_0,y_0)}}  \chi (t)=d$.
\end{enumerate}
\end{cor}

We will denote by 
$X_{(x_0,y_0)}$ the set of all characters in the algebra 
${\cal R}_{(x_0,y_0)}$ 
determined by nets. Analogously, we denote by 
$X_{(x_0,y_0)}^i$ the characters determined by the corresponding 
condition $i$ (where, $i=1,2,3,4,5$) of Theorem \ref{Th2}.

\begin{thm} \label{Th4}
$\forall f(x,y)\in {\cal R}_{(x_0,y_0)} \quad 
\forall d\in D_{(x_0,y_0)}^f \quad
\exists \chi \in X^1_{(x_0,y_0)}  \chi (t)=d$.
\end{thm}
\bproof
If $f(x,y) \in {\cal R}_{(x_0,y_0)}$ and
$d \in D^f_{(x_0,y_0)}$ then by Corollary 3 of Theorem \ref{Th2} we can 
conclude that there exists 
$\chi \in X_{(x_0,y_0)}$ such that 
$\chi (t)=d$.

If $\chi \notin X^1_{(x_0,y_0)}$ then we show that for 
$i=2,3,4,5$ by the description of $\chi$ in terms of Theorem \ref{Th1} 
(or Theorem \ref{Th1'}), we can choose 
$\chi_f \in  X^1_{(x_0,y_0)}$ and $\chi_f(f)=\chi(f)=d$.

It suffices to consider the case where the character is described by 
conditions of Theorem 1 for $x_0\leq x$. Let
$f(x,y)=\frac{p(x,y)}{q(x,y)}$. Suppose that
$\chi \in X^2_{(x_0,y_0)}$ is determined by the curve with the equation 
$y=y_0+\Phi_N (x-x_0)+(x-x_0)^e$ where $e\notin \mathbf{Q}$. then
$$p(x, y_0 +\Phi_N(x-x_0)+(x-x_0)^e)=a(x-x_0)^{u+ne}+0((x-x_0)^{u+ne})$$ and
$$q(x,y_0+\Phi_N(x-x_0)^e)=b(x-x_0)^{v+me}+0((x-x_0)^{v+me}).$$
Now choose $e_{N+1} \in {\cal Q}$ sufficiently close to $e$ such that:
$e_N < e_{N+1}$, 
$$p(x,y_0+\Phi_N(x-x_0)+(x-x_0)^{e_{N+1}})=
p(x_0,y_0)+a(x-x_0)^{u+ne_{N+1}}+0((x-x_0)^{u+ne_{N+1}})$$
$$q(x,y_0+\Phi_N(x-x_0)+(x-x_0)^{e_{N+1}})=
q(x_0,y_0)+b(x-x_0)^{v+me_{N+1}}+0((x-x_0)^{v+me_{N+1}}.$$
If $\chi_f$ is defined by 
$y-y_0+\Phi_N(x-x_0)+(x-x_0)^{e^1_{N+1}}$ then
$\chi_f \in X^1_{(x_0,y_0)}$ and $\chi_f(f)=d$.

In the cases of $\chi_f \in X^3_{(x_0,y_0)}$ and 
 $\chi_f \in X^4_{(x_0,y_0)}$, the proof is as in the case of
$\chi_f \in X^2_{(x_0,y_0)}$. Finally, if
$\chi_f \in X^5_{(x_0,y_0)}$ then $\chi (f)=\lim_{N\to{\infty}} \chi_N(f)$,
 where $\chi_N \in X_{(x_0,y_0)}$. It follows from Lemma 5 that
for sufficiently big $N=N_f$ holds
$\chi (f) = \chi_{N_f}(t)$.
\hfill $\Box$
\eproof
\vskip 0.4cm
\noindent {\bf Corollary}
Let $f(x,y)$ be a rational function. Then 
$\lim_{(x,y)\to (x_0,y_0)} f(x,y)=a\in {\cal R}$ exists iff
there exist the limits
$$\lim_{x\to x_0} f(x,y_0+\Phi_N(x-x_0))=a, \quad 
\lim_{y\to y_0} f(x_0+\Phi_N(y-y_0),y)=a$$
for arbitrary $\Phi_N(t) \in ({\cal I})$.

\bproof
 Indeed, if $f\in {\cal R}_{(x_0,y_0)} $ then our assertion is a 
particular case of Theorem \ref{Th4}. If $f(x,y)$ is unbounded in every nbd of
$(x_0,y_0)$ then repeating the arguments 
in the proof of Theorem \ref{Th4} for
$d=\infty$, it is easy to check that there exists $\Phi_N(t)$ such that
one of the following facts hold
$$\lim_{x\to x_0, x_0 <x} f(x, y_0+\Phi_N(x-x_0))=\infty$$
$$\lim_{x\to x_0, x <x_0} f(x, y_0+\Phi_N(x-x_0))=\infty$$
$$\lim_{y\to y_0, y_0 <y} f(x_0+\Phi_N(y-y_0),y)=\infty$$
$$\lim_{y\to y_0, y <y_0} f(x_0+\Phi_N(y-y_0),y)=\infty$$
\hfill $\Box$
\eproof
\vskip 1cm

\section{Algebras  $\bar {\cal R}_G$ and $\bar {\cal R}_{\bar G}$} 

Now we apply our results for description of characters to the case 
of certain Banach algebras which contain dense subalgebras 
of rational functions.
For these purposes we recall some known facts 
from the theory of Banach algebras \cite{[2]}.
Let $G$ be a connected open subset of $\mathbf{R}^2$ having compact 
closure $\bar G$. Denote by $B(\bar G)$ the Banach algebra of all
real bounded functions on $\bar G$ endowed with $sup$-norm.

By Gelfand-Naimark Theorem every closed subalgebra $A$ of $B(\bar G)$ 
is isomorphic with the Banach algebra $C(X^A)$ of all continuous 
functions on the compactum $X^A$, where $X^A$ is the set of all 
nontrivial characters of $A$ endowed with the weak topology 
$\sigma (X^A, A)$.

The desired isomorphism between $A$ and $C(X^A)$ is just the 
Gelfand Transform $\bar f (\chi)=\chi(f)$, where 
$f\in A, \bar f \in C(X^A)$ and $\chi \in X^A$.  

If a closed subalgebra $A$ separates the points of $B(\bar G)$ then 
identifying the character $\chi_{(x_0,y_0)} \in X^A$ (defined by 
$\chi_{(x_0,y_0)}(t)=f(x_0,y_0)$ ) with the point $(x_0,y_0)$, we may 
suppose that 
$\bar G$ is a dense subset of $X^A$. For instance, $\bar G$ is dense in 
$X^A$ if $C(\bar G) \subset A$.  

If $C(\bar G) \subset A$ then 
$X^A= \cup_{(x,y) \in \bar G} \chi^A_{(x,y)}$ where
$$\chi^A_{(x,y)}=
\{\chi \in X^A | \forall f \in C(\bar G) \quad  \chi(f)=f(x,y)\}.$$
Clearly $\chi^A_{(x,y)}$ is compact for every $(x,y) \in \bar G$.
Since $\bar G$ is dense in $X^A$ then 
$$\forall \chi \in X^A \quad \exists \{(x_{\gamma},
y_{\gamma})_{\gamma \in \Gamma} \quad \lim_{\gamma 
\in \Gamma}(x_{\gamma},y_{\gamma})=\chi$$
Moreover, if $\chi \in X^A_{(x,y)}$ then 
$\lim_{\gamma 
\in \Gamma}(x_{\gamma},y_{\gamma})=(x,y)$ in $\mathbf{R}^2$.

Note also that 
$D^{(x,y)}_f =\bar f (X^A_{(x,y)})$ for every $f \in A$.    

In the sequel we consider in $B(\bar G)$
the following two subalgebras:
$\bar {\cal R}_{\bar G}$ - the closure of a subalgebra consisting of all 
rational functions that are bounded on $\bar G$; and also 
$\bar {\cal R}_G$ -the closure of the subalgebra consisting of all rational
functions that are bounded on $G$ and continuous on the boundary of 
$G$. 
Obviously, 
$C(\bar G) \subset \bar {\cal R}_G \subset \bar {\cal R}_{\bar G}$. 

Let $X^{\bar G}$ ($X^G$) be the set of all characters of the algebra 
$\bar {\cal R}_{\bar G}$ (respectively, $\bar {\cal R}_G$) 
both endowed with the weak topology.

\begin{prop} \label{pr1}
The spaces $X^G$ and $X^{\bar G}$ are connected. Moreover, 
the sets $X^G_{(x,y)}$ and $X^{\bar G}_{(x,y)}$ are connected for every
 $(x,y) \in \bar G$. 
\end{prop}
   
This proposition is a consequence of Lemma \ref{lem2}.

\begin{prop} \label{pr2}
$\bar R_G$ {\it is not} separable and 
$X^G$ {\it is not} metrizable.
\end{prop}
\bproof
Assuming the contrary, let $\bar R_G$ be separable. Since the set of 
all rational functions is dense in 
$\bar R_G$, by \label{lem 1} we obtain that there are countably many 
points at which a function from $\bar R_G$ can be discontinuous. On the 
other hand, clearly every point may be a point of discontinuity for suitable 
function from $\bar R_G$. This contradiction proves the first assertion. 
Now, the Banach algebra $C(X^G)$ is not separable being topologically 
isomorphic to the Banach algebra $\bar R_G$. 
Therefore $X^G$ is not metrizable.
\hfill $\Box$
\eproof

\vskip 0.4cm    
\noindent {\bf Corollary}
$\bar R_{\bar G}$ {\it is not} separable and 
$X^{\bar G}$ {\it is not} metrizable.
\vskip 0.4cm
 
Note that Theorem \ref{Th2} provides a description of $X^G$ and $X^{\bar G}$.
Note also that $X_(x,y)\subset X^G_{(x,y)}=X^{\bar G}_{(x,y)}$ 
for any $(x,y) \in G$.
If $(x,y) \in \partial G$ ($\partial G$ means the boundary of $G$), 
then $\tilde X (x,y) \subset X^{\bar G}_{(x,y)}$, where 
$\tilde X (x,y)$ is the set of all characters in $X_{(x,y)}$, that can be 
determined by nets all members of which are in $\bar G$.

\begin{prop} \label {pr3}
$X^G=\cup_{(x,y) \in G} X_{(x,y)} \cup \partial G$ and 
$X^{\bar G}=\cup_{(x,y) \in G} X_{(x,y)} \cup 
(\cup_{(x,y) \in \partial G} {\tilde X}_{(x,y)})$.
\end{prop}    
\bproof
As we already mentioned, 
$X^G=\cup_{(x,y) \in \bar G} X^G_{(x,y)}$. Consider 
$(x_0.y_0) \in G$ and $\chi \in X^G_{(x_0,y_0)}$. 
There exists a net $\{(x_{\gamma},y_{\gamma})\}_{\gamma \in \Gamma}$ in $G$ 
which converges to $\chi$ in weak topology. Then this net determines a 
character from $X_{(x_,y_0)}$. Indeed, if not, as in the proof of 
Theorem \ref{Th2}, there exist natural $n_0, m_0$ such that the limit 
$\lim_{\gamma \in \Gamma} f(x_{\gamma},y_{\gamma})$ does not exist, 
where 
$$f(x,y)=\frac{(x-x_0)^{n_0}(y-y_0)^{m_0}}{(x-x_0)^{2n_0}+(y-y_0)^{2m_0}}$$
This is a contradiction because 
$f(x,y) \in {\cal R}_G$. Therefore, 
$\chi \in X_{(x_0,y_0)}$ and 
$X^G_{(x_0,y_0)}=X_{(x_0,y_0)}$, as desired. 
Analogous proof is valid for the second statement.
\hfill $\Box$
\eproof

\vskip 0.4cm
\noindent {\bf Corollary 1}
$X_{(x_0,y_0)}=X^{{\cal R}_{(x_0,y_0)}}$ for every $(x_0,y_0) \in {\bf R^2}$.
\vskip 0.4cm

\noindent {\bf Corollary 2}
$\bar G$ is sequentially dense in $X^G$, that is, for every $\chi \in X^G$ 
there exists a sequence which converges to $\chi$ and consist of points 
from $\bar G$.
\vskip 0.4cm 

It is actually a reformulation of Corollary 2 of Theorem \ref{Th2}.

\vskip 0.4cm
\noindent {\bf Corollary 3}
$X^1_{(x,y)}$ is dense in $X^{\bar G}_{(x,y)}=X_{(x,y)}$ for all 
$(x,y) \in G$ and 
${\tilde X}^1_{(x,y)}$ is dense in $X^{\bar G}_{(x,y)}=\tilde X_{(x,y)}$ 
for all $(x,y) \in \partial G$.

\vskip 0.4cm
 
This assertion follows from Theorem \ref{Th4}.

\vskip 0.4cm
\noindent {\bf Corollary 4}
Let $f(x,y) \in {\bar {\cal R}}_{\bar G}$ and $(x_0,y_0) \in \bar G$. 
Then 
$\lim_{(x,y) \to (x_0,y_0)} f(x,y) =a$ exists iff 
$\lim_{(x,y) \in L \to (x_0,y_0)} f(x,y) =a$ exists for every curve 
$L\subset \bar G$ that are described in the assertion 1 of Theorem 
\ref{Th1} 
and the assertion 1 of Theorem \ref{Th1'}.
\vskip 0.4cm

The verification is easy by the previous corollary and Corollary 1 
of Theorem \ref{Th2}.

\vskip 0.5cm

Let $f(x,y) \in B(\bar G)$ satisfies the following conditions:

\begin{enumerate}
\item
$f(x,y)$ is continuous on the boundary of $G$;
\item
\begin{itemize}
\item [a)]
For every $(x_0,y_0)\in G$ and a curve $L$ that is described 
in one of the assertions (1),(2),(3),(4) of Theorem 
\ref{Th1} or Theorem \ref{Th1'}  there 
exists the limit $$\lim_{(x,y)\in L \to (x_0,y_0)} f(x,y);$$
\item [b)]
For every $(x_0,y_0) \in G$ and a sequence of the curves 
$\{L_N\}_{N \in {\bf N}}$ that are described in Theorem \ref{Th1} 
or in Theorem \ref{Th1'}, 
there exists 
$$\lim_{N \to \infty} \lim_{(x,y) \in {L_N} \to (x_0.y_0)} f(x,y).$$
\end{itemize}
\end{enumerate}

\vskip 0.4cm

The set of all such functions clearly is a closed subalgebra of the 
Banach algebra $B(\bar G)$. Define this subalgebra by $L(\bar G)$. 
Obviously, $\bar {\cal R}_G \subset L(\bar G)$. However, 
$L(\bar G) \neq \bar {\cal R}_G$ as it follows from the following 
known counterexample  
$$f(x,y)=\frac{e^{-\frac{1}{(x-x_0)^2}}(y-y_0)}{e^{-\frac{2}{(x-x_0)^2}}+
(y-y_0)^2}$$
(see \cite{[1]}).

\vskip 0.4cm
\noindent {\bf Corollary 5}
If $f(x,y) \in \bar {\cal R}_G$ then:
\begin{enumerate}
\item 
The set of discontinuity points of $f(x,y)$ is at most countable;
\item
$f(x,y) \in L(\bar G)$;
\item
For every $(x_0,y_0) \in G$ and $d \in D^f_{(x_0,y_0)}$ there exists a 
character $\chi \in X^G_{(x_0,y_0)}$ 
such that 
$\chi (f)=d$.
\end{enumerate}

\vskip 0.5cm

We state here two natural questions inspired by Corollaries 4 and 5 of 
Proposition 3.

\begin{enumerate}
\item
How can be generalized Theorem 3 for 
$\bar {\cal R}_G$ or $\bar {\cal R}_{\bar G}$ ?
\item
Is it true that the conditions of Corollary 5 
are also {\it sufficient} 
in order to ensure that 
$f(x,y) \in \bar {\cal R}_{\bar G}$ ?
\end{enumerate} 

\vskip 0.3cm
Note that the main results of the present paper can be generalized 
for functions of $n$ arguments with $n > 2$ and for arbitrary 
points of discontinuity of rational functions.

\vskip 1cm
Finally the author thanks to E. Shustin (Tel-Aviv University) 
for several stimulating conversations and important suggestions 
that have significant influence on 
results of the present paper. We thank also to M. Megrelishvili 
(Bar-Ilan University) for his support.

\vskip 2cm
address:  \hskip 0.2cm Neot Golda str. 606/14, Netanya 42345, Israel
\vskip 0.2cm
tel: (972) 09-8356839
\vskip 0.3cm
\tt email:megereli@macs.biu.ac.il

\end{document}